\documentclass[arma]{svjour}
\usepackage{amsmath,amssymb,epsfig}

\numberwithin{equation}{section}

\newcommand{\ga}{\alpha}

\newcommand{\gev}{\varepsilon}

\newcommand{\bbr}{\mathbb{R}}

\newcommand{\bbn}{\mathbb{N}}

\newcommand{\bbe}{\mathbb{E}}

\newcommand{\opc}{\operatorname{C}}

\newcommand{\opl}{\operatorname{L}}

\newcommand{\oph}{\operatorname{H}}

\newcommand{\opw}{\operatorname{W}}

\newcommand{\opid}{\operatorname{id}}

\newcommand{\dive}{\operatorname{div}}

\newcommand{\dhr}{\mathrel{\lhook\joinrel\relbar\kern-.8ex\joinrel\lhook\joinrel\rightarrow}}

\begin{document}

\title{A new well-posed nonlocal Perona-Malik type equation}
\author{Patrick Guidotti}
\titlerunning{A new Perona-Malik type equation}

\maketitle

\begin{abstract}
A modification of the Perona-Malik equation is proposed for which the local nonlinear diffusion 
term is replaced by a nonlocal term of slightly reduced ``strength''. The new equation is globally 
well-posed (in spaces of classical regularity) and possesses desirable properties from the 
perspective of image processing. It admits characteristic functions as (formally linearly ``stable'') 
stationary solutions and can therefore be reliably employed for denoising keeping blurring in 
check. Its numerical implementation is stable, enhances some of the features of Perona-Malik, and 
avoids problems known to affect the latter.
\end{abstract}

\section{Introduction}
The Perona-Malik equation \cite{PM90}
\begin{equation}\label{pm}
  u_t-\dive\bigl(\frac{1}{1+|\nabla u|^2}\nabla u\bigr)=0
\end{equation}
on a domain $\Omega\overset{o}{\subset}\bbr ^n$ is usually complemented with homogeneous Neumann boundary 
conditions. The nonlinear diffusion coefficient
$$
 a(u)=\frac{1}{1+|\nabla u|^2}
$$ 
is sometimes replaced by $a(u)=e^{-|\nabla u|^2}$. Perona and Malik originally introduced \eqref{pm} 
in the context of image processing with the aim of denoising a given image $u_0$ while at the same 
time controlling blurring. The latter is unavoidable if a linear diffusion equation is used instead. 
The nonlinear diffusion determined by $a(u)$ significantly impedes diffusion in directions of steep 
gradients which correspond to sharp edges in the image. If one considers a large gradient across a 
surface, then the nonlinear diffusion coefficient leads to diffusion in tangential directions and 
backward diffusion in normal direction. This is the reason why models of this type are often referred 
to as anisotropic diffusion models. \\
Numerical approximations of \eqref{pm} have confirmed this prediction and have been successfully 
employed as a consequence. It has, however, been observed that discretizations of \eqref{pm} often lead to 
stair-casing in the numerical solution (see Figure \ref{eps}), or, more in general, to the preservation 
of gradients which do not represent any image feature but might simply be due to the presence of noise 
(cf. \cite{PM90,ALM92} and Figure \ref{kink_osc}). 
Equation \eqref{pm} poses many a challenge from the analytical point of view. It is an example of 
so-called forward-backward diffusion. In \cite{K97} the author shows that it is not well-posed. To 
circumvent this problem a variety of regularization techniques have been proposed. Spatial regularizations 
were the first to appear in \cite{CLMC92,ALM92}. The gain in analytical tractability is, however, offset by 
the introduction of blurring. Temporal and spatio-temporal regularizations have also been considered, 
see \cite{NS92,CEA98,CB01,B03,BC05,Ama06}. In \cite{Ama06} the author considers a purely temporal 
regularization which leads to time-delayed Perona-Malik equation which mimics the implicit linearization 
procedure common to all numerical discretizations 
of \eqref{pm}, that is
$$\begin{cases}
 u^{n+1}=u^n+\dive \bigl( a(u^n)\nabla u^{n+1}\bigr)\delta t&n\geq 0\, ,\\ 
 u^0=u_0 &n=0\, ,\end{cases}
$$
whereby the nonlinearity is always evaluated at the previous time step. The introduction of a time 
delay leads to a tractable analytical problem and does not seem to negatively impact the desired 
features of a corresponding numerical implementation. The author proves local existence of a regular 
solution for the delay equation.\\ 
Variational techniques have also been utilized to gain analytical understanding of \eqref{pm} 
in a one-dimensional context. The Perona-Malik equation can in fact be viewed as the gradient flow 
associated with the non-convex energy functional
$$
  E(u)=\int _0^1 \log\bigl( 1+|\nabla u|^2\bigr)\, ,\: u\in\oph ^1(\Omega)\, .
$$
Using concepts and techniques related to Young-measure solutions various results have been 
obtained. \cite{TTZ05,Z06} prove instability of certain solutions and the existence of infinitely 
many solutions of a certain (weak) type. \\
Equation \eqref{pm} was motivated by its potential applications to image processing. Whereas it 
poses significant and interesting mathematical challenges, its form is clearly not dictated by 
any physical principle. A modification is therefore proposed in this paper which is globally well-posed, 
and allows for special natural functions to be stationary solutions. The latter leads to a desirable 
dynamical behavior from the practical point of view. It also leads to stable pseudo-spectral 
discretizations with 
satisfactory properties from the perspective of image processing. It can also be viewed as a new 
regularization technique. It differs, however, substantially from other regularization techniques in that 
it actually makes piecewise constant function become stationary solutions and in that the degree of 
regularization can be more finely tuned. In particular it does not regularize at any given specific 
scale.\\
To simplify the discussion, analytical considerations will be restricted to the one dimensional 
setting. A representative numerical experiment, however, will be presented also in the natural two 
dimensional setting. Roughly speaking, the modification proposed consists in replacing the nonlinear term 
by
\begin{equation}\label{mpf}
  a_\varepsilon(u)=\frac{1}{1+[(-A)^{\frac{1-\gev}{2}}u]^2}  
\end{equation}
for $\varepsilon\in(0,\frac{1}{2})$ and where $A$ is the Neumann Laplacian. The 
fractional powers appearing in \eqref{mpf} can be defined in various manners depending on the 
choice of function spaces in which the equation is considered.
Since one is interested in denoising while preserving sharp edges, it is clear that $a_\gev$ should 
provide the very same benefits for the kind of large gradients one is interested in preserving. 
The associated equation 
\begin{equation}\label{mpm1}
  u_t-\bigl(\frac{1}{1+[(-A)^{1-\varepsilon}u]^2}u_x\bigr) _x=0
\end{equation}
has, however, the advantage of being quasi-linear and of avoiding backward diffusion altogether. 
The analysis and the numerical experiments performed in this paper therefore also show that backward 
diffusion is not an essential ingredient in order to preserve sharp edges. For analytical reasons 
the ideas just presented will be performed on an equivalent formulation of \eqref{pm} leading to 
an equation of the flavor of \eqref{mpf}.
\section{The Problem}
In its one dimensional formulation, the Perona-Malik equation reads
\begin{equation}\label{1pm}
 \begin{cases} u_t-\bigl(\frac{1}{1+u_x^2}u_x\bigr)_x=0 &\text{ in }(0,1)\text{ for }t>0\, ,\\
 u_x=0&\text{ at }x=0,1\text{ for }t>0\, ,\\
 u=u_0&\text{ at }t=0\text{ in }\, [0,1].
\end{cases}
\end{equation}
Introducing $v(x)=\int_{0}^{x}\!u(y)\,dy$ as a new unknown, it is easily seen that it satisfies
\begin{equation}\label{1pmi}
 \begin{cases}
  v_t-\frac{1}{1+v_{xx}^2}v_{xx}=0&\text{ in }(0,1)\text{ for }t>0\, ,\\
  v(t,0)=0\, ,\: v(t,1)=\int_{0}^{1}\!u_0(y)\,dy&\text{ for }t>0\, ,\\
  v(0,\cdot)=\int_{0}^{\cdot}\!u_0(y)\,dy&\text{ in }[0,1]\, ,
 \end{cases}
\end{equation}
since it is readily verified that the average of $u_0$ is preserved over time using \eqref{1pm}. 
By further defining 
$w(x)=v(x)-x\int_{0}^{1}\!u_0(y)\,dy\, ,\: x\in[0,1]\, ,$ it follows that
\begin{equation}\label{1pmf}
  \begin{cases}
  w_t-\frac{1}{1+w_{xx}^2}w_{xx}=0&\text{ in }(0,1)\text{ for }t>0\, ,\\
  w(t,0)=0=w(t,1)&\text{ for }t>0\, ,\\
  w(0,\cdot)=\int_{0}^{\cdot}\!u_0(y)\,dy-(\cdot)\int_{0}^{1}\!u_0(y)\,dy&\text{ in }[0,1]\, .
  \end{cases}
\end{equation}
Conversely, any solution of \eqref{1pmi} leads to a solution of \eqref{1pm} by setting 
$$
 u(x)=v_x(x)\, ,\: x\in[0,1]\, .
$$
The evolution equation is clearly satisfied whereas for the boundary conditions the 
following argument is needed. The function 
$$
 \tilde v(x)=\int_{0}^{x}\!u(y)\,dy\, ,\: x\in[0,1]\, ,
$$ 
satisfies
$$
  \tilde v_t-\frac{1}{1+\tilde v _{xx}^2}\tilde v_{xx}=-\frac{u_x(0)}{1+u_x^2(0)}\, .
$$
Since $\tilde v=v$ by $v(0)=0$ and since $v$ satisfies \eqref{1pmi}, it follows that 
$\frac{u_x(0)}{1+u_x^2(0)}=0$ and therefore $u_x(0)=0$. Now
$$
 \int_{0}^{1}\!u(t,y)\,dy=v(t,1)-v(t,0)=\int_{0}^{1}\!u_0(y)\,dy
$$
implies that
$$
 0=\int_{0}^{1}\!u_t(y)\,dy=\frac{u_x(1)}{1+u_x^2(1)}-\frac{u_x(0)}{1+u_x^2(0)}
$$
and therefore that $u_x(1)=0$, too. Notice that if $u$ is a piecewise constant function, then 
$v$ is a continuous piecewise affine function. Thus piecewise affine functions play the same 
role for \eqref{1pmf} as piecewise constant functions do for \eqref{1pm}.\\
Equation \eqref{1pmf} is fully nonlinear and clearly presents the same analytical difficulties 
as the original Perona-Malik equation. The following modification is proposed
\begin{equation}\label{mpm}
  \begin{cases}
  u_t-a_\varepsilon(u)Au=0&\text{ in }(0,1)\text{ for }t>0\, ,\\
  u(0,\cdot)=u_0&\text{ in }[0,1]\, ,
  \end{cases}  
\end{equation}
for 
\begin{equation}\label{mpff}
 a_\varepsilon(u)=\bigl(1+[(-A)^{1-\varepsilon}u]^2\bigr)^{-1}\, ,\: 
 \varepsilon\in(0,1/2)\, ,
\end{equation} 
and where $A$ now denotes the Dirichlet Laplacian. This is perfectly analogous to \eqref{mpf}. 
Notice that the intensity of the nonlinearity is barely reduced. It is therefore to be expected 
that solutions of \eqref{mpm} with $\varepsilon>0$ behave similarly to the solutions of \eqref{mpm} 
with $\varepsilon=0$ for which the original \eqref{pm} is recovered, at least at the discrete level. 
At the continuous level \eqref{1pm} is ill-posed and such a claim is not very meaningful.
It is therefore likely that \eqref{mpm} can provide a viable practical denoising tool while preventing 
blurring.
\begin{remark}\label{rem:cons}
It is easily verified that
$$
 \frac{d}{dt}\int_{0}^{1}\!v\,dx=0\text{ and }\frac{1}{2}\frac{d}{dt}\int_{0}^{1}\!|v|^2\,dx=
 -\int_{0}^{1}\!\frac{|v_x|^2}{1+[A^{1-\varepsilon}u]^2}\,dx
$$
for $v=u_x$ and a smooth solution $u$ of \eqref{mpm}.
\end{remark}

\section{Global Existence of Smooth Solutions}
It turns out that \eqref{mpm} can be shown to be well-posed in the classical sense as soon as 
$\varepsilon>0$. In order to prove this, some notation needs to be introduced. For $\alpha\in[0,1)$ 
let 
$$
 \opc ^\alpha _0(0,1):=\begin{cases}
 \{ u\in\opc ^\alpha\bigl( [0,1],\bbr\bigr)\, |\, u(j)=0\, ,\: j=0,1\} &\alpha\in(0,1)\, ,\\
 \opc\bigl([0,1],\bbr\bigr)& \alpha=0\, ,\end{cases}
$$
where $\opc ^\alpha\bigl( [0,1],\bbr\bigr)$ is for $\alpha>0$ the standard space of H\"older continuous 
functions endowed with norm
$$
 \| \cdot\| _\alpha:=\| \cdot\| _\infty +[\cdot]_\alpha\, .
$$
The semi-norm $[\cdot ]_\alpha$ is given by
$$
  [u]_\alpha:=\sup _{x\neq y}\frac{|u(x)-u(y)|}{|x-y|^\alpha}\, ,\: u\in\opc ^\alpha\bigl( [0,1],
  \bbr\bigr)
$$
if $\alpha>0$. Given $\alpha\in[0,1)$, the operator $A_\ga$ is given by
$$
  \begin{cases} 
    \operatorname{dom}{A_\alpha}=\opc ^{2+\alpha}_0(0,1):=\{ u\in\opc ^{2+\alpha}(0,1)\, |\, 
    u,u_{xx}\in\opc ^\alpha _0(0,1)\} &\text{ for }\alpha>0\\
    \operatorname{dom}{A_0}=\opc ^2_0(0,1):=\{ u\in\opc ^2(0,1)\, |\, u(j)=0\, ,\: j=0,1\}&
    \text{ for }\alpha=0
  \end{cases}
$$
and
$$
  A_\alpha u=u_{xx}\, ,\: u\in\operatorname{dom}(A_\alpha)\, .
$$
Once $\alpha\in[0,1)$ is fixed, the following simplified notation will be used
$$
  A:=A_\alpha\, ,\: E_0:=\opc ^\alpha _0(0,1)\, ,\: E_1=\opc ^{2+\alpha}_0(0,1)\, .
$$
It follows from \cite[Corollary 3.1.21, Corollary 3.1.32, Theorem 3.1.34]{Lun95} that $A$ is a sectorial 
operator and therefore generates an analytic semi-group $\{e^{tA}\, |\, t\geq 0\}$ on $E_0$. Then its 
fractional powers $(-A)^\gamma$, $\gamma\in(0,1)$, can be defined as the inverses of
$$
 (-A)^{-\gamma} :=\frac{1}{\Gamma(\gamma)}\int_{0}^{\infty}\!t^{\gamma -1}e^{tA}\,dt
$$
defined on their range, that is, by
$$
 \operatorname{dom}\bigl((-A)^\gamma\bigr)=\operatorname{R}\bigl((-A)^{-\gamma}\bigr)^{-1}\, ,\:
 (-A)^\gamma =\bigl((-A)^{-\gamma}\bigr)^{-1}\, .
$$
The domain of $(-A)^\gamma$ might not be an interpolation space between $E_1$ and $E_0$ but it 
satisfies the interpolation inequality
\begin{equation}\label{int-ineq}
 \| (-A)^\gamma x\| _{E_0}\leq c\| x\| _{E_0}^{1-\gamma}\| Ax\| _{E_0}^\gamma\, ,\: x\in
 \operatorname{dom}(A)\, .
\end{equation}
It can, however, be sandwiched between known interpolation spaces
\begin{equation}\label{sandwich}
 (E_0,E_1)_{\gamma,1}\hookrightarrow D\bigl((-A)^\gamma\bigr)\hookrightarrow(E_0,E_1)_{\gamma,\infty}
\end{equation}
where, for $\gamma\in(0,1)$ and $p\in[1,\infty]$, $(E_0,E_1)_{\gamma,p}$ is the standard real 
interpolation functor. Observe that for $\alpha=0$ it can proved that
$$
 (E_0,E_1)_{\gamma,\infty}=\opc ^{2\gamma}_0(0,1)\, ,\: \gamma\in(0,1)\setminus\{\frac{1}{2}\}\, .
$$
H\"older spaces of the time variable will also be useful. For $\beta\in(0,1)$, let
\begin{equation}\label{mrfs}
  \bbe _0:=\opc ^\beta\big([0,T],E_0\bigr)\, ,\\ \bbe _1:=\opc ^{1+\beta}\bigl([0,T],E_0\bigr)
  \cap\opc ^{\beta}\bigl([0,T],E_1\bigr)\, .
\end{equation}
If $u\in\bbe _1$, then
$$
 \frac{1}{1+[(-A)^{1-\varepsilon}u(t)]^2}A
$$
is sectorial for every  fixed $t$. By H\"older maximal regularity (see \cite[Proposition 6.1.3]{Lun95}) 
it follows that any solution of
$$\begin{cases}
 v_t-\frac{1}{1+[(-A)^{1-\varepsilon}u]^2}Av=0 &\text{ for }t>0\\
 v(0)=u_0&\end{cases}
$$
satisfies $u\in\bbe _1$ for any $u_0\in E_1$ such that $Au_0\in (E_0,E_1)_{\beta,\infty}$. In order 
to apply the mentioned theorem, one needs also to take into account the additional facts that
$$
 [t\mapsto A(t):=\frac{1}{1+[(-A)^{1-\varepsilon}u(t)]^2}A]\in\opc ^\beta\bigl([0,T],\mathcal{L}
 (E_1,E_0)\bigr)
$$
and 
$$
 \operatorname{dom}\bigl( A(t)\bigr)=\operatorname{dom}\bigl(A(0)\bigr)\, ,\: t\in[0,T]\, .
$$
It is now possible to prove the following result.
\begin{theorem}\label{thm:ex}
Let $u_0\in E_1$ such that $Au_0\in (E_0,E_1)_{\beta,\infty}$. Then there exists $T^+(u_0)$ and a unique 
$u:\bigl[0,T^+(u_0)\bigr)\to E_1$ such that $u\big |_{[0,T]}$ is a solution of \eqref{mpm} on 
the time interval $[0,T]$ for any $0<T<T^+(u_0)$. Furthermore
$$
 u\in\opc ^{1+\beta}\bigl( [0,T^+(u_0)),E_0\bigr)\cap\opc ^{\beta}\bigl( [0,T^+(u_0)),E_1\bigr)
$$
\end{theorem}
\begin{proof}
Let $v\in\bbe _1$ and define $\Phi(v)$ to be the solution of
$$
 u_t-\frac{1}{1+[(-A)^{1-\varepsilon}v]^2}Au=0\, ,\: u(0)=u_0\, .
$$
Define $r=2\| u_0-\Phi(u_0)\| _{\bbe _1}$ and let
$$
 B_r:=\{ v\in\bbe _1\, |\, \| v-u_0\| _{\bbe _1}\leq r\}\, .
$$
It follows that
\begin{multline*}
  \| \Phi(v)-u_0\| _{\bbe _1}\leq\| \Phi(v)-\Phi(u_0)\| _{\bbe _1}+\| \Phi(u_0)-u_0\| _{\bbe _1}=\\
  \big\|\bigl[ \partial _t-a_\varepsilon(v)A\bigr]^{-1}(0,u_0)-\bigl[ \partial _t-a_\varepsilon(u_0)A
  \bigr]^{-1}(0,u_0)\big\| _{\bbe _1}+\frac{r}{2}\\ 
  =\Big\|\bigl[ \partial _t-a_\varepsilon(u_0)A\bigr]^{-1}\Big\{ 
  1-\bigl[ \partial _t -a_\varepsilon(u_0)A\bigr]\bigl[ \partial _t-a_\varepsilon(v)A\bigr]^{-1}
 \Big\}(0,u_0)\Big\| _{\bbe _1}+\frac{r}{2}\\\leq r
\end{multline*}
if $T>0$ is chosen small enough since $v\in B_r$ and $v(0)=u_0$, so that $v$ is uniformly close 
to $u_0$. To see this, consider
\begin{multline}\label{eq1}
 \|[a_\varepsilon(v)-a_\varepsilon(u_0)]A\| _{\mathcal{L}(\bbe _1,\bbe _0)}\leq
 \| a_\varepsilon(v)-a_\varepsilon(u_0)\| _{\mathcal{L}(\bbe _0)}\\\leq\|(-A)^{1-\varepsilon}(v-u_0)
 \|_{\bbe _0} 
\end{multline}
which follows from
$$
 \| \frac{1}{1+w^2(t)}-\frac{1}{1+w^2(s)}\| _{E_0}\leq \| w(s)-w(t)\|_{E_0}
$$
by setting $w=(-A)^{1-\varepsilon}v$ since $v(0)=u_0$ and $a_\varepsilon(u_0)$ are independent of 
the time variable.
To further estimate \eqref{eq1}, one can use the interpolation inequality \eqref{int-ineq} and 
\begin{equation}\label{eq2}
 v(t)-v(s)=\int _s^t\dot v(\tau)\, d\tau\text{ and } \dot v\in\opl _\infty E_0\, 
\end{equation}
to obtain
\begin{multline*}
  \| (-A)^{1-\varepsilon}[v(t)-u_0]\| _{E_0}\leq c\|v(t)-u_0\| _{E_1}^{1-\varepsilon}\| 
  v(t)-u_0\| _{E_0}^\varepsilon\\\leq crt^{\beta(1-\varepsilon)+\varepsilon}\, ,\: s,t\leq T\, ,
\end{multline*}
and
\begin{multline*}
  \| (-A)^{1-\varepsilon}[v(t)-v(s)]\| _{E_0}\leq c\| v(t)-v(s)\| _{E_1}^{1-\varepsilon}
 \| v(t)-v(s)\| ^\varepsilon _{E_0}\\\leq cr|t-s|^{\beta(1-\varepsilon)+\varepsilon}\, ,\: s,t\leq T\, .
\end{multline*}
It follows that
$$
 \bigl[(-A)^{1-\varepsilon}(v-u_0)\bigr] _{\beta}\leq cr T^{(1-\beta)\varepsilon}\, .
$$
and finally that
$$
 \|[a_\varepsilon(v)-a_\varepsilon(u_0)]A\| _{\mathcal{L}(\bbe _1,\bbe _0)}\leq crT^{(1-\beta)
 \varepsilon}
$$
Thus the claimed inequality follows by a simple Neumann series argument made possible by choosing 
$T$ small enough. Next observe that $\Phi(v_1)-\Phi(v_2)$ solves
\begin{equation*}\begin{cases}
 [\Phi(v_1)-\Phi(v_2)] _t=a_\varepsilon(v_1)A[\Phi(v_1)-\Phi(v_2)]+\bigl( a_\varepsilon(v_1)-
 a_\varepsilon(v_2)\bigr)A\Phi(v_2)\\ [\Phi(v_1)-\Phi(v_2)](0)=0
\end{cases}
\end{equation*}
It follows that
$$
 \Phi(v_1)-\Phi(v_2)=\bigl[ \partial _t-a_\varepsilon(v_1)A\bigr]^{-1}\Big\{ 0,\bigl[ a_\varepsilon(v_1)-
 a_\varepsilon(v_2)\bigr]A\Phi(v_2)\Big\}\, .
$$
Thus
\begin{multline*}
  \|\Phi(v_1)-\Phi(v_2)\| _{\bbe _1}\leq c(r)\|[a_\varepsilon(v_1)-a_\varepsilon(v_2)]A\Phi(v_2)\| _{E_0}
  \\\leq c(r)\| A^{1-\varepsilon}[v_1-v_2]\| _{\bbe _0}\| \Phi(v_2)\| _{\bbe _1}
\end{multline*}
Now, since $v_1-v_2\in\bbe _1$ it follows that
$$
 \| A^{1-\varepsilon}[v_1-v_2]\| _{\bbe _0}\leq c\| v_1-v_2\| _{\bbe _1} T^{\varepsilon(1-\beta)}
$$
by making use of the interpolation inequality \eqref{int-ineq} and of \eqref{eq2} in a perfectly 
similar way to the calculations used to obtain the self-map property.
Thus Banach fixed-point Theorem implies the existence of a unique solution 
$$
 u\in\bbe _1
$$ 
provided the time interval is chosen short enough. It follows from \cite[Proposition 6.1.3]{Lun95} that 
$$
  \dot u(t)\in (E_0,E_1)_{\beta,\infty}\, ,\: t\in (0,T]
$$ 
which gives $Au(T)\in(E_0,E_1)_{\beta,\infty}$. Thus, by using the same argument, the solution can be 
extended to a larger interval of existence. Repeating this argument indefinitely, a solution is obtained 
on a maximal interval of existence $\bigl[ 0,T^+(u_0)\bigr)$ with the stated properties.
\end{proof}
Next classical Sobolev spaces are needed
$$
 \opw ^2_p(0,1)=\big\{ u\in\opl _p(0,1)\, \big |\, \partial^j u\in\opl _p(0,1)\text{ for }
 0\leq j\leq 2\big\}
$$
for $p\in(0,\infty)$ where the derivatives have of course to be understood in the distributional 
sense. The simple embedding
\begin{equation}\label{rkineq}
  \opw ^s_p(0,1)\hookrightarrow\opc ^{s-1/p}(0,1)\, ,\: sp>1\, ,
\end{equation}
will be very useful in the proof of the next Lemma.
\begin{lemma}\label{lem:glob-ex}
If it can be shown that $u\in\opl _\infty\bigl(\bigl[ 0,T^+(u_0)\bigr),\opw ^2_p(0,1)\bigr)$ for all 
$p\in(1,\infty)$, then the solution exists globally in time.
\end{lemma}
\begin{proof}
Observe that
$$
 \dot u=a_\varepsilon(u)\triangle u\in\opl _\infty\bigl(\bigl[ 0,T^+(u_0)\bigr),\opl _p(0,1)\bigr)
$$
follows from the assumption and the form of $a_\varepsilon$. For $s,t\in\bigl[ 0,T^+(u_0)\bigr]$, one has
\begin{multline*}
 \| A^{1-\varepsilon}\bigl( u(t)-u(s)\bigr)\| _{\opc ^{\tilde{\rho}}}\\\leq c\,
 \| A^{1-\varepsilon}\bigl( u(t)-u(s)\bigr)\| _{\opw ^{2\rho}_p}\leq c\,
 \| A^{1-\varepsilon+\rho}\bigl( u(t)-u(s)\bigr)\| _{\opl _p}\\\leq
 c\, \| A\bigl( u(t)-u(s)\bigr)\| _{\opl _p}^{1-\varepsilon+\rho}\|\int_{s}^{t}\!\dot u(\tau)\,d\tau
 \| _{\opl _p}^{\varepsilon-\rho}\leq c\, |t-s|^{\varepsilon-\rho}
\end{multline*}
as follows from \eqref{rkineq} and \eqref{int-ineq}. Hereby it needs to be assumed that
$$
 \tilde\rho =2\rho-\frac{1}{p}>0\text{ and }0<\rho<\varepsilon\, .
$$
This can always be achieved by choosing $p$ large enough and yields
\begin{equation}\label{reg}
 A^{1-\varepsilon}u\in\opc ^{\varepsilon-\rho}\bigl(\bigl[ 0,T^+(u_0)\bigr],\opc ^{\tilde{\rho}}(0,1)
 \bigr)\, .  
\end{equation}
Denote the function space in \eqref{reg} by $\bbe _0$. Let then $v$ be the solution of
$$
 \dot v -a_\varepsilon(u)Av=0\, ,\:v(0)=u_0
$$
on $\bigl[ 0,T^+(u_0)\bigr]$. It satisfies $v\in\bbe _1$ and one has
$$
 v\big | _{[0,T]}=u\big |_{[0,T]}\, ,\: T<T^+(u_0)
$$
by uniqueness. Since $v$ is uniformly (H\"older) continuous, so must be $u$. It can therefore be 
extended as a solution to $\bigl[ 0,T^+(u_0)\bigr]$ and as a consequence of Theorem \ref{thm:ex} beyond 
any finite time. The solution is therefore global. 
\end{proof}
\begin{lemma}\label{lem:smooth}
Let $u_0\in\opc^\infty\bigl([0,1]\bigr)$ satisfy compatibility conditions to all orders, that is, 
assume that $u_0\in\operatorname{dom}\bigl((-A)^k\bigr)$, $k\in\bbn$. Then the solution of \eqref{mpm} 
satisfies
$$
 u\in\opc ^{1+\beta}\Bigl(\bigl[ 0,T^+(u_0)\bigr),\opc ^\infty\bigl([0,1]\bigr)\Bigr)\, ,
$$
for $\beta<\varepsilon$.
\end{lemma}
\begin{proof}
Let $u$ be the solution of \eqref{mpm}. Then 
$$
 (-A)^{1-\varepsilon}u\in\opc ^\beta\Bigl(\bigl[ 0,T^+(u_0)\bigr),\opc ^{\alpha+2\beta}_0(0,1)
\Bigr)
$$ 
since $\beta<\varepsilon$ and $u$ is therefore also a solution of \eqref{pm} in $\bbe _1$ for 
$\alpha+\beta$. By repeating this boot-strapping argument indefinitely one obtains the claim.
\end{proof}
\begin{remark}\label{rem:timereg}
Time regularity can also be improved by similar boot-strapping arguments but it will play no role in 
this paper. 
\end{remark}
It is now possible to prove the following global regularity results.
\begin{theorem}\label{thm:mthm}
The maximal solution of \eqref{mpm} exists globally.
\end{theorem}
\begin{proof}
According to Lemma \ref{lem:glob-ex} it is enough to show 
$$
 \| Au(t)\| _{\opw ^2_p}\leq c\, ,\: t\in\bigl[ 0,T^+(u_0)\bigr)\text{ for }p\in(0,\infty)\, .
$$
By Lemma \ref{lem:smooth} it follows from \eqref{pm} that
$$\begin{cases}
 A\dot u-A\bigl( a_\varepsilon(u)Au\bigr)=0\, ,\\
 Au(0)=u_0\, .
\end{cases}
$$
Thus $v:=Au$ satisfies
\begin{equation}\label{Aeq}\begin{cases}
 \dot v-A\bigl( a_\varepsilon(u)v\bigr)=0\, ,\\
 v(0)=u_0\, .
\end{cases}
\end{equation}
Observe that $a_\gev(u(t))A\cdot$ is, for any fixed time $t$, formally adjoint to 
the operator $A\bigl[ a_\varepsilon\bigl(u(t)\bigr)\cdot\bigr]$. 
The evolution operator $U_{\mathbb{A}}(t,\tau)$ generated by $\mathbb{A}=a_\varepsilon(u)A$ on $E_0$ 
satisfies
$$
 \| U_{\mathbb{A}}(t,\tau)u_0\| _{\opl _p(0,1)}\leq \| u_0\| _{\opl _p(0,1)}\, .
$$
This follows from the Trotter-Kato product formula (see \cite{Pa83}) and the fact that the operator
$$
 c(x)A
$$
generates a contraction semigroup on $\opl _p(0,1)$ if ellipticity is assumed (cf. \cite{A07}).  
Let $T<T^+(u_0)$ and consider the family of generators $\mathbb{B}=\mathbb{A}(T-t)$, $t\in[0,T]$. 
Then an easy computation based on
$$
 \frac{d}{d\tau}U_{\mathbb{A}}(t,\tau)u_0=-U_{\mathbb{A}}(t,\tau)\mathbb{A}(\tau)u_0\, ,\:
 u_0\in\operatorname{dom}{\mathbb{A}(0)}\, ,
$$
reveals that
$$
 U^*_{\mathbb{B}}(T,T-t)v_0\, ,\: t\in[0,T]\, ,
$$
satisfies \eqref{Aeq} if $U^*_{\mathbb{B}}$ on $\opl _{p'}(0,1)$ is taken to be the evolution operator 
dual to $U_{\mathbb{B}}$. Now, since
$$
 \| U^*(t,\tau)\| _{\mathcal{L}(\opl _{p'})}=\| U(t,\tau)\| _{\mathcal{L}(\opl _p)}\leq 1
$$
it follows that
$$
 \| v(t)\| _{\opl _{p'}}\leq \| v_0\| _{\opl _{p'}}=\| Au_0\| _{\opl_{p'}}
$$
for any $p\in(1,\infty)$ and therefore 
$$
 \| Au(t)\| _{\opl _p}\leq c\,\| u\| _{\opw ^2_p}\leq c\,\| Au_0\| _{\opl _p}
$$
for all $p\in(1,\infty)$.
The proof is thus complete since $T<T^+(u_0)$ is arbitrary and the embedding constant $c$ does not 
depend on $T$.
\end{proof}
\section{Stationary Solutions}
Among constant functions, $u\equiv 0$ is clearly  a stationary solution of \eqref{mpm}. It 
corresponds to a constant solution with value $\int _0^1 u_0(x)\, dx$ of the original Perona-Malik 
equation with initial value $u_0$. It turns out, however, that other more desirable 
functions are equilibria of \eqref{mpm} as well. Even though the next result remains valid for \eqref{mpm}, 
it is formulated in the context of the modified Perona-Malik equation with periodic boundary conditions 
instead of homogeneous Dirichlet conditions. This allows the simplify a little the argument.
\begin{proposition}\label{sationary}
Let $u$ be a continuous $1$-periodic piecewise affine function. Then
$$
 \frac{1}{1+[(-A)^{1-\varepsilon}u]^2}Au\equiv 0\, ,
$$
for $A$ taken to be the realization of $\partial_{xx}$ in the space of Radon measures $M_\pi$. 
The subscript $\pi$ indicates that the periodic case is considered. 
\end{proposition}
\begin{proof}
If $u$ has the properties stated in the proposition, it follows that
$$
 \partial_{xx}u=\sum _{j=1}^nc_j\delta _{x_j}
$$ 
for some constants $c_j\in\bbr$ and locations $x_1,\dots,x_n\in[0,1)$. The claim would therefore 
follow if it can proved that
$$
 \frac{1}{1+[(-A)^{1-\varepsilon}u]^2} 
$$
is a continuous function vanishing at $x_j\, ,\: j=1,\dots,n$.
It is clear that it is enough to consider the case of a function with a single kink. Furthermore, since 
periodicity is assumed, Poisson's summation formula implies that the singularity in 
$(-A)^{1-\varepsilon}u$ in [0,1) is the same as the singularity at the same location for 
$(-A)^{1-\varepsilon}u$ periodically extended to real line. By making use of standard localization 
techniques, it is therefore sufficient to 
consider a function of the real line with a single kink in the origin and $A=\partial_{xx}$ on the 
full line. Choose the prototype $u(x)=|x|$. Since
$$
  (-A)^{1-\varepsilon}u=(-A)^{-\varepsilon}\delta
$$
it follows that $\mathcal{F}\bigl((-A)^{1-\varepsilon}u\bigr)=\frac{1}{|\xi |^{2\varepsilon}}$ and 
therefore that
$$
 (-A)^{1-\varepsilon}u=c_\varepsilon\frac{1}{|x|^{1-2\varepsilon}}\, .
$$
for $C_\varepsilon=\sqrt{\frac{2}{\pi}}\Gamma(1-2\varepsilon)\sin(\pi\varepsilon)$. 
Remember that it has been assumed initially that $\varepsilon\in(0,1/2)$ to make the nonlinearity 
strong enough. It follows that, in the general case, $(-A)^{1-\varepsilon}u$ has integrable 
singularities at $x_j$, $j=1,\dots,n$, and is otherwise smooth by hypo-ellipticity. This shows 
that
$$
 \frac{1}{1+[(-A)^{1-\varepsilon}u]^2}
$$
is a continuous function which vanishes exactly at $x_j$, $j=1,\dots,n$. It is therefore a well-defined 
multiplier for $\sum _{j=1}^nc_j\delta _{x_j}$ and one has
$$
 \frac{1}{1+[(-A)^{1-\varepsilon}u]^2}Au=0
$$
as claimed.
\end{proof}
\begin{remark}\label{rem:norealstationary}
It should be observed that it is not possible (to the best of our knowledge) to make sense of 
\eqref{pm} as generating a flow on $M_\pi$ (or on any of the standard Besov spaces containing Dirac 
distributions) and so the stationary solutions of the proposition are only formally such. 
A short computation would even show that they are formally linearly 
stable. This might help explain the fact that, in numerical calculations of \eqref{mpm}, it is 
observed that solution starting close to such a stationary solution tend to stay in its vicinity 
for a long time before being driven to the trivial steady state. Even more is actually true. Numerical 
experiments (see Figure \ref{regev}) show that smooth solutions tend to become piecewise affine at first 
but are eventually completely smoothed out. 
\end{remark}
\begin{remark}\label{rem:linearization}
Since the weakened nonlinearity \eqref{mpff} does not give rise to forward-backward diffusion for the 
linearization, global well-posedness can be proved. Its linearization, however, will be closer and closer 
to the Perona-Malik linearization as $\gev$ decreases. As a consequence more and more of its eigenvalues 
will become negative when linearizing in the presence of large gradients. This is supporting evidence 
that the numerical benefits of Perona-Malik should not go lost in the modified \eqref{mpm}.  
\end{remark}
\begin{remark}\label{rem:asym-behavior}
The new equation is globally well-posed for smooth initial conditions. It is, however, not immediately 
clear what their large time asymptotic behavior should be. Numerical experiments seem to suggest that 
any solution will eventually converge to a trivial steady-state. This behavior might, however, be due 
to numerical diffusion. In the next section a reason will be given which seems to exclude the latter 
possibility.
\end{remark}
\begin{remark}\label{rem:stability}
It is also observed that the original Perona-Malik equations do not admit piece-wise constant 
functions as formal stationary points unless an ad-hoc concept of generalized stationary solution 
is introduced (see \cite{K97}) since the multiplication of distributions is in general not well-defined. 
Its discretization, however, seems to exhibit ``almost stability'' of such solutions in the sense 
explained here. This might be an indication that the Perona-Malik equation is regularized by 
discretization and could explain the success of its implementations in spite of the ill-posedness of 
its analytical counterpart.
\end{remark}
\section{Numerical implementation and experiments}
\begin{figure}
\begin{center}
\includegraphics[scale=0.2]{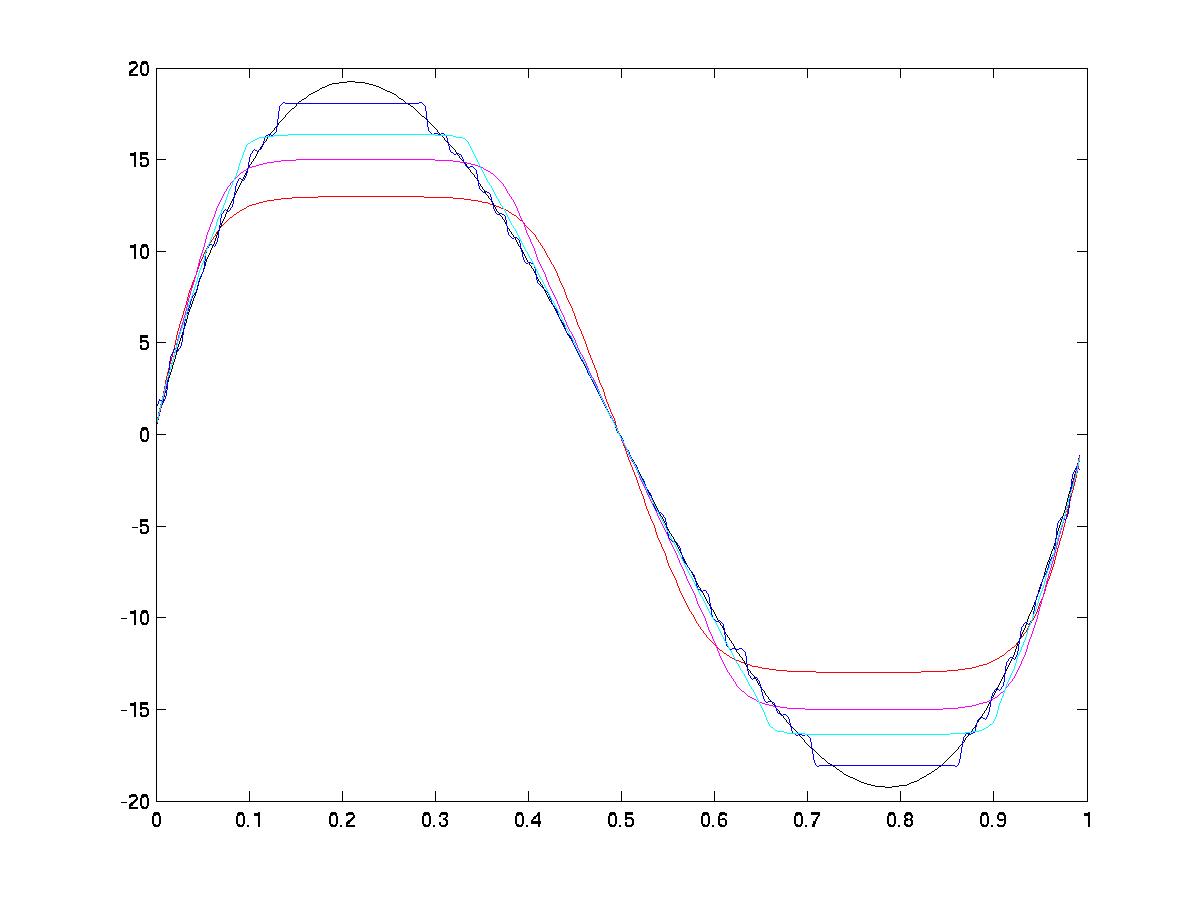}
\caption{The gradient of the solution for different values of $\varepsilon$.}\label{eps}
\end{center}
\end{figure}

\begin{figure}
\begin{center}
\includegraphics[scale=0.5]{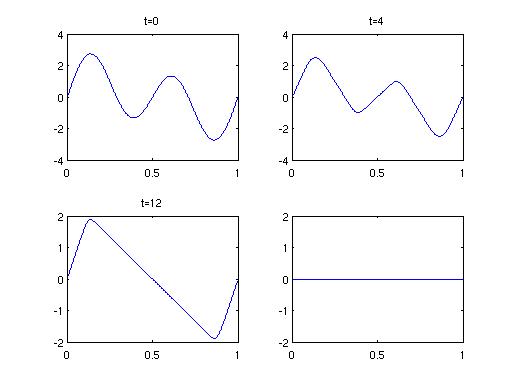}
\caption{Tendency to evolve towards or close to piece-wise affine functions.}\label{regev}
\end{center}
\end{figure}

\begin{figure}
\begin{center}
\includegraphics[scale=0.2]{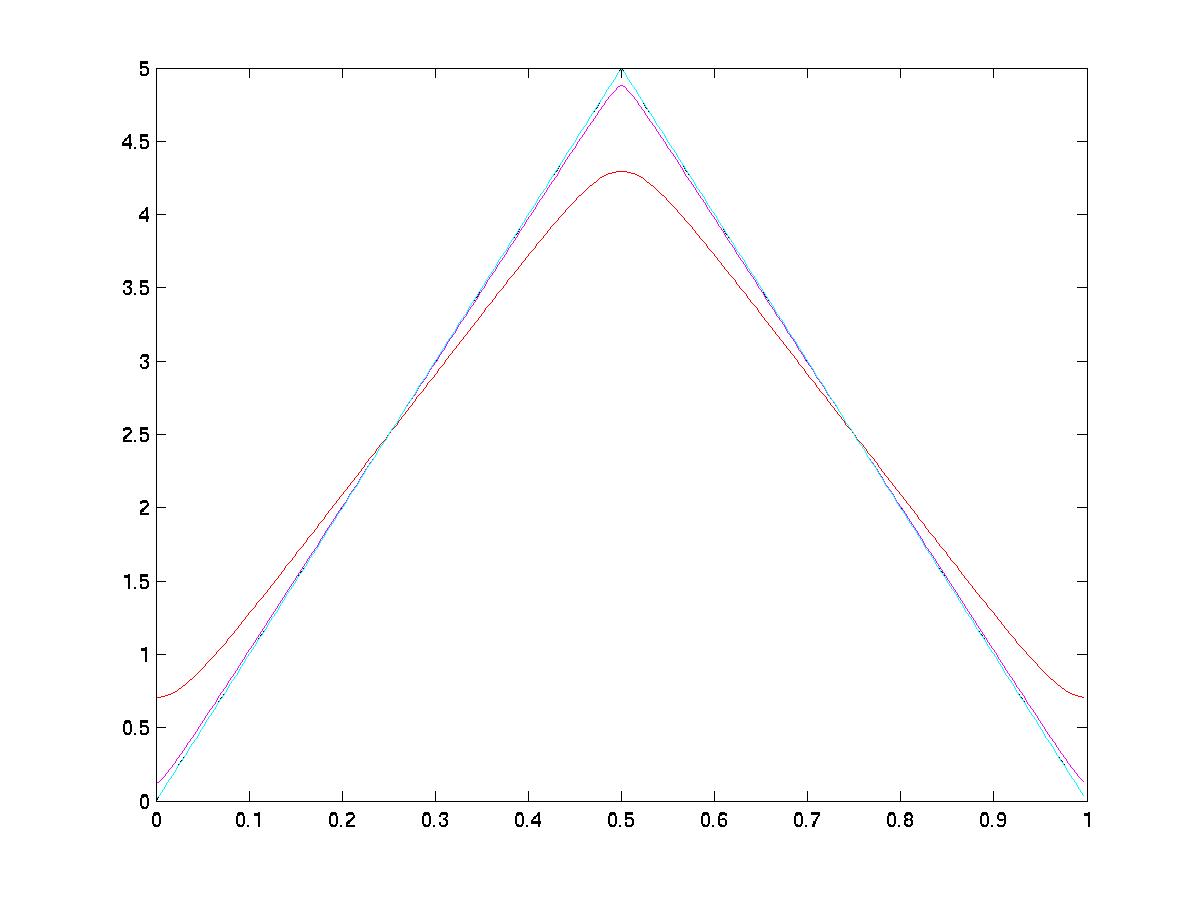}
\caption{Behavior close to formally stationary solutions.}\label{kink} 
\end{center}
\end{figure}

\begin{figure}
\begin{center}
\includegraphics[scale=0.5]{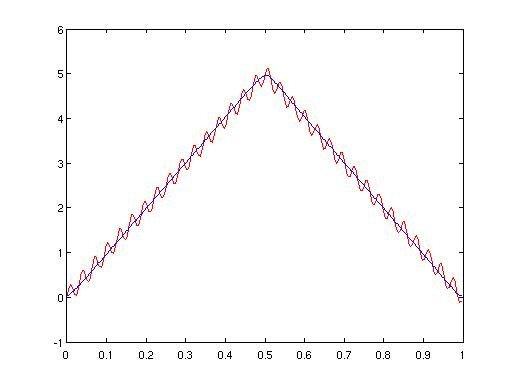}
\caption{The effect of \eqref{mpm} on high frequency low contrast oscillations. The initial condition 
(in red) is rapidly evolved to the almost piecewise affine function (in blue).}\label{kink_osc}
\end{center}
\end{figure}

\begin{figure}
\begin{center}
\includegraphics[scale=0.6]{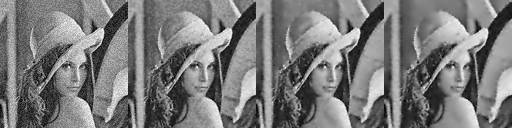}
\caption{The denoising effect obtained for the 2D version of \eqref{mpm} with 
         $\varepsilon=0.6$ and Neumann boundary conditions. The initial condition is Lenna's image 
         corrupted with about 15\% salt and pepper noise.}\label{teaser}
\end{center}
\end{figure}

Next a numerical discretization of \eqref{mpm} is derived in a periodic context and used to perform 
numerical experiments intended to illustrate and demonstrate the claimed improvements on the classical 
Perona-Malik equation. The periodic Laplacian $A$ is discretized spectrally by means of the fast Fourier 
transform $F_n$
\begin{equation}\label{plaplacian}
  A_n=F_n^{-1}\Lambda _nF_n
\end{equation}
where $n=2^m$ denotes the number of grid points used and
$$
 D_n=4\pi ^2\operatorname{diag}((-\frac{n}{2}+1)^2,(-\frac{n}{2}+2)^2,\dots, 
 0,1,\dots,\frac{n^2}{4})\, .
$$
The time variable is discretized by forward semi-implicit Euler so that
\begin{equation}\label{time-stepping}
  u^{k+1}=\bigl[\opid _n-\frac{h_t}{1+(A_n^{1-\varepsilon}u^k)^2}A_n  \bigr]^{-1}u^k
\end{equation}
where $u^k$ is the spatial $n$-vector at time $k\, h_t$ for the time step $h_t>0$. Observe 
that, setting $\varepsilon=0$, the classical Perona-Malik equation is recovered.\\
The following experiments are presented here. First the evolution of the function
$$
  u_0(x):=100x^2(1-x^2)\, ,\: x\in[0,1]\, .
$$
is considered. Figure \ref{eps} depicts the first derivative of the function $u^k$ after $100$ 
time steps of size $h_t=0.06$ for $m=8$. The blue curve corresponds to the Perona-Malik solutions 
and the stair-casing 
phenomenon is apparent. The red, magenta and cyan curves correspond to solution of the modified 
equation for $\varepsilon=0.3,0.2,0.1$, respectively. They clearly benefit from all effects of the 
Perona-Malik equation without leading to stair-casing. Some decrease in contrast is the only price 
paid. Notice that the black curve depicts the first derivative of the initial value. Comparison with 
the other curves reveals the de-blurring effect of the equation.
It has been proved that continuous piecewise affine functions are steady-states for the modified 
Perona-Malik equation. There is, however, no mechanism by which a smooth solution can develop a 
singularity in finite time. The solution therefore feels the presence of these nontrivial steady-states 
but never leaves the regime where diffusion, slowly but surely, drives it to a trivial steady-state. 
This is exemplified in Figure \ref{regev} where various stages of the evolution of a smooth solution are 
with initial condition
$$
 u_0(x)=\sin(2\pi x)+2\sin(4\pi x)\, ,\: x\in[0,1]\, .
$$
is depicted. Clearly the convergence towards the trivial steady state might be due to numerical 
diffusion (see Remark \ref{rem:asym-behavior}). However, Remark \ref{rem:cons} provides an indirect 
way to test the numerics. If follows from Remark \ref{rem:cons} that
$$
 \int_{0}^{1}\!u_x^2\,dx=\int_{0}^{1}\!(u_0)_x^2\,dx-2\int_{0}^{t}\!\int_{0}^{1}\!
 \frac{u_{xx}^2}{1+[A^{1-\varepsilon}u]^2}\,dx\,d\tau\, .
$$
This conservation relation can be tracked for any smooth solution and, in the particular case considered, 
the relative numerical deviation observed between the left and right-hand-side of it amounts to a mere 
$0.2$\% for a fully converged solution.
Scheme \eqref{time-stepping} also preserves the piecewise affine structure of initial values for 
large integration times. The solution with initial condition
$$
 u_0(x)=5-|10\,x-5|\, ,\: x\in(0,1)\, ,
$$
is computed for $100$ time steps with $h_t=0.06$ and $m=8$. The solutions are plotted in Figure \ref{kink}. 
The color coding is as in the previous experiment. In spite of the fact that $u_0$ is a steady-state 
for any choice of $\varepsilon\in(0,\frac{1}{2})$, numerical dissipation is stronger for larger 
$\varepsilon$ as the relative strength of the non-linearity decreases. The initial value and the solution 
to $\varepsilon=0.1$ are indistinguishable in the plot. The differ by $0.2\%$ in the supremum norm. 
Continuous piecewise linear functions are not steady-states for the original Perona-Malik equation since 
the non-linearity is not well-defined for such functions. In spite of this its numerical counterpart 
delivers results comparable to those for small positive $\varepsilon$.
Figure \ref{kink_osc} shows clearly one of the claimed enhancements of \eqref{mpm} over the original 
Perona-Malik equations. The initial condition (in red)
$$
 u_0(x)=5-|10\,x-5|+0.2\sin(64\pi x)\, ,\: x\in(0,1)\, ,
$$
is evolved to time $t=2$ (in blue) with $\varepsilon=0.3$. The new equation can manifestly 
differentiate between high low contrast gradients and high contrast gradients, which are remarkably 
well preserved. Observe that the initial oscillatory condition would be left virtually unchanged by 
the original Perona-Malik equation for the same and longer time ranges.
Even though this paper is concerned with a one-dimensional modification of Perona-Malik, equation 
\eqref{mpm} can be considered in a natural two dimensional setting. The qualities of the its 
one dimensional counterpart considered here do carry over to that case. Figure \ref{teaser} 
shows the evolution of a noisy test image every pixel of which has been corrupted by about 15\% noise 
in the gray-scale. For other tests and details of the two-dimensional implementation we refer to 
\cite{GL07}.
\section{Conclusions}
A modified Perona-Malik equations has been proposed by strength reduction in the non-linearity via 
a non-local term. Global well-posedness of the new equation and the fact that piecewise constant 
functions are steady-states for it are its main analytical advantages. From the practical point 
of view, the new equation delivers enhanced benefits as compared to Perona-Malik and suppresses known 
shortcomings associated to it. It provides an easy, effective, and stable tool for image de-noising. 
These claims are corroborated by numerical experiments.

\bibliographystyle{plain}
\bibliography{../../lite}

\begin{thebibliography}{10}

\bibitem{ALM92}
L.~Alvarez, P.-L. Lions, and J.-M. Morel.
\newblock Image selective smoothing and edge-detection by non-linear diffusion.
  ii.
\newblock {\em SIAM J. Numer. Anal.}, 29(3):845--866, 1992.

\bibitem{Ama06}
H.~Amann.
\newblock Time-{D}elayed {P}erona-{M}alik {P}roblems.
\newblock {\em Acta Math. Univ. Comenianae}, LXXVI:1--24, 2006.

\bibitem{A07}
W.~Arendt.
\newblock Heat kernels.
\newblock {\em Internet Seminar 2005/2006,\\
  https:tulka.mathematik.uni-ulm.de/2005/lectures/ internetseminar.ps}.

\bibitem{B03}
A.~Belahmidi.
\newblock {\em {E}quations aux d\'eriv\'ees partielles appliqu\'ees \`a la
  restoration et \`a l'aggrandissiment des images.}
\newblock Ph.D. Thesis. Universit\'e Paris-Dauphine, Paris, 2003.

\bibitem{BC05}
A.~Belahmidi and A.~Chambolle.
\newblock Time-delay regularization of anisotropic diffusion and image
  processing.
\newblock {\em {M2AN Math. Model. Numer. Anal.}}, 39(2):231--251, 2005.

\bibitem{CLMC92}
F.~Catt\'e, P.-L. Lions, J.-M. Morel, and T.~Coll.
\newblock Image selective smoothing and edge-detection by non-linear diffusion.
\newblock {\em SIAM J. Numer. Anal.}, 29(1):182--193, 1992.

\bibitem{CB01}
Y.~Chen and P.~Bose.
\newblock On the incorporation of time-delay regularization into
  curvature-based diffusion.
\newblock {\em J. Math. Imaging Vision}, 14(2):149--164, 2001.

\bibitem{CEA98}
G.~Cottet and M.~El Ayyadi.
\newblock A volterra type model for image processing.
\newblock {\em IEEE Trans. Image Processing}, 7:292--303, 1998.

\bibitem{GL07}
P.~Guidotti and J.~Lambers.
\newblock {A Well-posed Nonlinear Nonlocal Diffusion for Noise Reduction}.
\newblock {\em Submitted}.

\bibitem{K97}
S.~Kichenassamy.
\newblock The {P}erona-{M}alik paradox.
\newblock {\em SIAM J. Appl. Math.}, 57(5):1328--1342, 1997.

\bibitem{Lun95}
A.~Lunardi.
\newblock {\em {A}nalytic {S}emigroups and {O}ptimal {R}egularity in
  {P}arabolic {P}roblems}.
\newblock Birk\-h{\"a}user, {B}asel, 1995.

\bibitem{NS92}
M.~Nitzberg and T.~Shiota.
\newblock Nonlinear image filtering with edge and corner enhancement.
\newblock {\em IEEE Trans. Pattern Anal. and Machine Intelligence},
  14:826--833, 1992.

\bibitem{Pa83}
A.~Pazy.
\newblock {\em {S}emigroups of {L}inear {O}perators and {A}pplication to
  {P}artial {D}ifferential {E}quations}.
\newblock Springer Verlag, New York, 1983.

\bibitem{PM90}
P.~Perona and J.~Malik.
\newblock Scale-space and edge detection using anistotropic diffusion.
\newblock {\em IEEE Transactions Pattern Anal. Machine Intelligence},
  12:161--192, 1990.

\bibitem{TTZ05}
S.~Taheri, Q.~Tang, and K.~Zhang.
\newblock Young measure solutions and instability of the one-dimensional
  {P}erona-{M}alik equation.
\newblock {\em J. Math. Anal. Appl.}, 308:467--490, 2005.

\bibitem{Z06}
K.~Zhang.
\newblock Existence of infinitely many solutions for the one-dimensional
  {P}erona-{M}alik model.
\newblock {\em Calc. Var.}, 26(2):171--199, 2006.

\end{thebibliography}

\address{103 Multipurpose Science and Technology Building\\ Department of Mathematics\\ 
         University of California\\ Irvine, CA 92697-3875 USA\\ email: {gpatrick@math.uci.edu}}

\end{document}